\documentclass[a4paper,twoside,10pt]{article}

\usepackage{amssymb}
\usepackage{amsmath}
\usepackage{amscd}
\usepackage[matrix]{xypic}
\usepackage{cite}

\newtheorem{theorem}{Theorem}

\newtheorem{definition}[theorem]{Definition}

\newtheorem{proposition}[theorem]{Proposition}

\newcommand\K{\mathbb{K}}

\begin{document}

\title{A class of nonassociative algebras and cogebras}
\author{Michel Goze\thanks{Email address: M.Goze@uha.fr}}
\date{}
\maketitle

\begin{center}
Laboratoire de Math\'ematiques et Applications, \\
Universit\'{e} de Haute Alsace,\\
4, rue des Fr\`eres Lumi\`ere, F-68093 Mulhouse Cedex, France.
\end{center}

\bigskip

 
\medskip

{\small \ \ \ \    {\bf MS Classification numbers}: 16Axx}

\bigskip

\abstract{We present classes of nonassociative algebras whose associator satisfies 
invariance conditions given by the action of the 3 order symmetric group. Amongst these 
algebras we find the wellknown Vinberg algebras, the Pre-Lie algebras, 
the Lie-admissible algebras and the 3-power associative algebras.}

\bigskip

\begin{center}
{\bf I. $\Sigma_3$-ASSOCIATIVE ALGEBRAS}
\end{center}

\medskip

\noindent {\bf A. $\Sigma _3$-invariant spaces}

\medskip

Let $\Sigma_3$ be the $3$-order symmetric group and $\mathbb{K}\, [\Sigma_3]$ the associated group algebra,
where $\mathbb{K}$ is a field of characteristic zero.
 We denote by  $\tau_{ij}$ be the transposition echanging $i$ and $j$, $c_1=(1,2,3)$ and $ c_2={c_1}^2$ the two $3$-cycles of $\Sigma_3.$
Every $v \in \mathbb{K}\,[\Sigma_3]$ decomposes as follows:
$$v=a_1 id + a_2 \tau _{12}  + a_3 \tau _{13} + a_4 \tau _{23}+a_5 c_1+a_6 c_2$$
or simply
$$v=\sum_{\sigma \in \Sigma_3} a_{\sigma} \sigma
$$
where $a_{\sigma} \in \mathbb{K}.$

Consider the natural right action of $\Sigma_3$ on $\mathbb{K}\,[\Sigma_3]$
$$
\begin{array}{ccc}
\Sigma_3 \times \mathbb{K}\, [ \Sigma_3 ] & \rightarrow & \mathbb{K}\,[ \Sigma_3 ] \\
 (\sigma, \sum_i a_i \sigma_i ) & \mapsto & \sum_i a_i \sigma^{-1} \circ \sigma_i
\end{array}
$$
For every vector $v \in \mathbb{K}\,[\Sigma_3]$
we denote by $\mathcal{O}(v)$ the corresponding orbit of $v.$ Let $F_v=\mathbb{K}\,(\mathcal{O}(v))$ be the linear subspace
of $\mathbb{K}\,[\Sigma_3]$ generated by $\mathcal{O}(v).$ It is an invariant subspace of $\mathbb{K}\,[\Sigma_3].$
Therefore, using Mashke's theorem, we deduce that it is a direct sum of irreducible invariant subspaces.

\medskip

\noindent {\bf B. $\Sigma _3$-associative algebras}

\medskip

Let $(\mathcal{A},\mu)$ be a $\mathbb{K}$-algebra with multiplication $\mu$.
We denote by $A_{\mu}$ the associator of $\mu$ that is
$$A_{\mu}= \mu \circ (\mu \otimes Id - Id \otimes \mu).
$$
Every $\sigma \in \Sigma_3$ defines a linear map denoted by $\Phi_{\sigma}$ given by
$$
\begin{array}{rccc}
\Phi_{\sigma} : & \mathcal{A}^{\otimes ^3} & \rightarrow & \mathcal{A}^{\otimes ^3} \\
& x_1 \otimes x_2 \otimes x_3 & \rightarrow & x_{\sigma^{-1}(1)} \otimes x_{\sigma^{-1}(2)} \otimes x_{\sigma^{-1}(3)}.
\end{array}
$$
If $v=\sum_{\sigma \in \Sigma_3} a_{\sigma} \sigma $ is a vector of $ \mathbb{K}\,[\Sigma_3]$ 
we define the endomorphism $\Phi_v$
of $\mathcal{A}^{\otimes ^3}$ by taking:
$$\Phi_v=\sum_{\sigma \in \Sigma_3} a_{\sigma} \Phi_{\sigma}.$$

\begin{definition}
An algebra $(\mathcal{A},\mu)$ is a $\Sigma _3$-associative algebra if there exists $v \in \mathbb{K}\, [\Sigma_3]$ such that
$$ A_{\mu} \circ \Phi_v =0.
$$

\end{definition}

\medskip

\noindent {\bf C. Lie-admissible algebras. $3$-power associative algebras}

\medskip

\begin{proposition}
Let $v$ be in $\mathbb{K}\,[\Sigma_3]$ such that  $dim \ F_v=1$. Then $v=\alpha V$ 
or $v= \alpha W$ with $\alpha \in \K$
where the vectors $V$ and $W$ are the following vectors :
\begin{eqnarray}
\label{1} V= Id - \tau _{12} -  \tau _{13} - \tau _{23}+ c_1 +c_2, \\
\label{2} W= Id + \tau _{12} +  \tau _{13} + \tau _{23}+ c_1 +c_2. 
\end{eqnarray}
\end{proposition}
The first case corresponds to the character of $\Sigma_3$ given by the signature, the second 
corresponds to the
trivial case.

\medskip  

\noindent Every algebra  $(\mathcal{A},\mu)$ whose associator satisfies
$$ A_{\mu} \circ \Phi_V =0
$$
is a Lie-admissible algebra. This means that the algebra  $(\mathcal{A},[,])$ whose product is given by the bracket $[x,y]=\mu(x,y)-\mu(y,x)$ is 
a Lie algebra. Likewise
an algebra  $(\mathcal{A},\mu)$ whose associator satisfies
$$ A_{\mu} \circ \Phi_W =0
$$
is $3$-power associative that is it satisfies $A_{\mu}(x,x,x)=0$ for every $x \in A.$

\bigskip

\begin{center}
{\bf II. $G_i$-ASSOCIATIVE ALGEBRAS}
\end{center}

\medskip

In this section we study the $\Sigma_3$-associative algebras corresponding to the subgroups of $\Sigma_3$.

\medskip

\noindent {\bf A. Notations}

\medskip

Let us consider the subgroups of $\Sigma_3$:
$$
\begin{array}{l}
G_{1}=\{Id\}, \\
G_{2}=\{Id,\tau _{12}\}, \\
G_{3}=\{Id,\tau_{23}\}, \\
G_{4}=\{Id,\tau _{13}\},  \\
G_{5}=\{Id,c_1,c_2\}=A_3 \  \mbox{\rm(the alternating group)} , \\
G_{6}=\sum_{3}
\end{array}
$$
We  keep these notations for the subgroups of $\Sigma _3$ in the following sections. 

\medskip

\noindent {\bf B. $G_i$-associative algebras}

\medskip
This notion was introduced by E.Remm in her thesis [10].

\begin{definition}
Let $G_i$ be a subgroup of ${\sum}_{3}.$ The algebra $(\mathcal{A},\mu )$ is $G_i$-associative if
\[
\sum_{\sigma \in G_{i}}(-1)^{\varepsilon (\sigma )}A_{\mu }\circ \Phi_{\sigma}
=0.
\]
\end{definition}
\begin{proposition}
Every $G_i$-associative algebra is a $\Sigma _3$-associative algebra.
\end{proposition}
{\it Proof.} Every subgroup $G_i$ of $\Sigma _3$ corresponds to an invariant linear space 
$F(v_i)$ generated by a single vector 
$v_i \in \Bbb{K}[\Sigma _3]$.
More precisely we consider $v_1=Id$, $v_2=\tau _{12}$, $v_3=\tau _{23}$, $v_4=\tau _{13}$, 
$v_5=Id+c_1+c_2$ and $v_6=V$ that we have defined in $(1)$. 
$\Box$
\begin{proposition}
Every $G_i$-associative algebra is a Lie-admissible algebra.
\end{proposition}
{\it Proof.} The vector $V$ belongs to the orbits ${\cal{O}}(v_i)$ for every $v_i$. Thus, if 
$\mu$ is a  $G_i$-associative product,
it also satisfies
$$A_{\mu }\circ \Phi_{V}=0
$$
and $\mu$ is a Lie-admissible multiplication. $\Box$

\medskip

\noindent We deduce the following type of Lie-admissible algebras :

\medskip

1. A $G_{1}$-associative algebra is an associative algebra.

2. A  $G_{2}$-associative algebra is a Vinberg algebra.
If $A$ is finite-dimensional, the associated Lie algebra is provided with an affine structure.

3. A  $G_{3}$-associative algebra is a  pre-Lie algebra.

4. If $(\mathcal{A},\mu) $ is $G_{4}$-associative then $\mu $ satisfies
\[
(X.Y).Z-X.(Y.Z)=(Z.Y).X-Z.(Y.X)
\]
with $X.Y=\mu (X,Y).$

5. If $(\mathcal{A},\mu) $ is $G_{5}$-associative then $\mu $ satisfies the generalized
Jacobi condition :
\[
(X.Y).Z+(Y.Z).X+(Z.X).Y=X.(Y.Z)+Y.(Z.X)+Z.(X.Y)
\]
with $X.Y=\mu (X,Y).$
Moreover if the product is antisymmetric, then it is a  Lie algebra bracket.

6. A  $G_{6}$-associative algebra is a Lie-admissible algebra.

\bigskip

\begin{center}
{\bf III. $G_i$-ASSOCIATIVE COGEBRAS}
\end{center}

\medskip

\noindent {\bf A. An axiomatic definition of the $G_i$-associative algebras}

\medskip

To define a $G_i$-cogebra it is easier to present an equivalent and axiomatic definition of the notion of $G_i$-associative
algebra.
\begin{definition}
A $G_i$-associative algebra is $(\mathcal{A},\mu,\eta,G_i)$ where $\mathcal{A}$ is a
vector space, $G_i$ a subgroup of $\Sigma_3$,
$\mu :\mathcal{A} \otimes \mathcal{A}  \longrightarrow \mathcal{A} $  and $\eta : \mathbb{K}
\rightarrow \mathcal{A}$ are linear maps satisfying the
following axioms :

\medskip

\noindent 1. ($G_i$-ass):  The square
$$
\begin{CD}
\mathcal{A} \otimes \mathcal{A}  \otimes \mathcal{A} @>{(\mu \otimes Id})_{G_i}>> \mathcal{A} \otimes \mathcal{A} \\
@V{(Id \otimes \mu)}_{G_i}VV  @V{\mu}VV \\
\mathcal{A} \otimes \mathcal{A}  @ >{\mu}>> \mathcal{A}
\end{CD}
$$
commutes, where ${(Id \otimes \mu)}_{G_i}$ is the linear mapping defined by:
$$
{(Id \otimes \mu)}_{G_i}=\sum_{\sigma \in G_i} (-1)^{\epsilon (\sigma )} (Id \otimes \mu) 
\circ \Phi_{\sigma}.
$$

\noindent 2. (Un) The following diagram is commutative :

$$
\xymatrix{ \mathbb{K} \otimes \mathcal{A} \ar[r]^{\eta \otimes id} \ar[dr]_{} &
\mathcal{A} \otimes \mathcal{A}  \ar[d]^{\mu}  & \mathcal{A} \otimes
\mathbb{K} \ar[l]_{id \otimes \eta} \ar[dl]{} \\  & \mathcal{A} & }
$$
\end{definition}

The axiom ($G_i$-ass) expresses  that the multiplication
$\mu$ is $G_i$-associative whereas the axiom (Un) means that the
element $\eta(1)$ of $\mathcal{A}$ is a left and right unit for $\mu$.

\begin{definition}
A morphism of $G_i$-associative algebras
$$f:(\mathcal{A},\mu,\eta,G_i)
\rightarrow (\mathcal{A}',\mu',\eta',G_i)$$
 is a linear map from $\mathcal{A}$ to $\mathcal{A}'$ such
that
$$ \mu' \circ (f \otimes f)=f \circ \mu \quad \mbox{and} \quad f \circ \eta=\eta'
$$
\end{definition}

\medskip

\noindent {\bf B. $G_i$-associative cogebras}

\medskip 

We want to dualize the previous diagrams to obtain the notions of corresponding cogebras.
Let $\Delta$ be a comultiplication on a vector space $C$:
$$\Delta : C \longrightarrow C \otimes C.$$
 We define the bilinear map
 $$G_i \circ (\Delta \otimes Id): C^{\otimes ^3}\longrightarrow C^{\otimes ^3}$$
 by
$$G_i \circ (\Delta \otimes Id)=\sum_{\sigma \in G_i} (-1)^{\epsilon (\sigma )} \Phi_{\sigma} \circ (\Delta \otimes Id).$$
\begin{definition}
A $G_i$-cogebra is a vector space $C$ provided with a comultiplication
$\Delta : C \longrightarrow C \otimes C$ and a counit $\epsilon: C
\rightarrow \mathbb{K}$ such that : 

\noindent 1. ($G_i$-ass co) The following square is commutative :

$$
\begin{CD}
C  @>{\Delta}>> C \otimes C \\
@V{\Delta}VV  @VV{G_i \circ (Id \otimes \Delta)}V \\
C \otimes C  @ >{G_i \circ (\Delta \otimes Id)}>> C \otimes C \otimes C.
\end{CD}
$$

\noindent 2. (Coun) The following diagram is commutative :

$$
\xymatrix{ \mathbb{K} \otimes C &
C \otimes C  \ar[l]_{\epsilon \otimes id} \ar[r]^{id \otimes \epsilon}  & C \otimes
\mathbb{K} \\  & C \ar[u]^{\Delta}  \ar[ur]{} \ar[ul]_{} & }
$$
\end{definition}

\begin{definition}
A morphism of $G_i$-associative cogebras 
$$f:(C,\Delta,\epsilon,G_i)
\rightarrow (C',\Delta',\epsilon',G_i)$$
 is a linear map from $C$ to $C'$ such
that
$$ (f \otimes f) \circ \Delta  =\Delta'  \circ f  \quad \mbox{and}
\quad \epsilon = \epsilon \circ f
$$
\end{definition}

\medskip

\noindent {\bf C. The dual space of a $G_i$-associative cogebra}

\medskip

For any natural number $n$ and any $\K$-vector spaces $E$ and $F$, we denote by
$$ \lambda _n : Hom(E,F)^{\otimes ^n}\longrightarrow Hom(E^{\otimes ^n},F^{\otimes ^n})$$
the natural embedding
$$\lambda _n(f_1\otimes ...\otimes f_n)(x_1\otimes ...\otimes x_n)=f_1(x_1)\otimes ...\otimes f_n(x_n).$$

\begin{proposition}
The dual space of a $G_i$-associative cogebra is a $G_i$-associative algebra.
\end{proposition}
{\it Proof}. Let $(C,\Delta)$ a $G_i$-associative cogebra. We consider the multiplication on the
dual vector space $C^*$ of $C$ defined by :
$$\mu= \Delta ^*\circ \lambda _2.$$
It provides  $C^{*}$
with a $G_i$-associative algebra structure. In fact we have
$$\mu (f_1 \otimes f_2)=\mu_{\K}\circ \lambda _2(f_1\otimes f_2)\circ \Delta \ \ \ (3) $$
for all $f_1 , f_2 \in C^{*}$ where $\mu_{\K}$ is the multiplication of $\K$. The equation $(3)$ 
becomes :
$$
\begin{array}{ll}
\mu\circ ( \mu \otimes Id)(f_1\otimes f_2\otimes f_3)& = \mu_{\K}\circ (\lambda _2(\mu(f_1\otimes f_2)\otimes f_3)\circ \Delta \\
\smallskip \\
& = \mu_{\K}\circ \lambda _2((\mu_{\K}\circ \lambda _2(f_1\otimes f_2)\circ \Delta)\otimes f_3)\circ \Delta \\
\vspace{0.001in} \\
&= \mu_{\K}\circ (\mu_{\K}\otimes Id)\circ \lambda _3(f_1\otimes f_2\otimes f_3)\circ (\Delta \otimes Id)\circ \Delta .
\end{array}
$$
The associator $A_{\mu }$ satisfies :
$$
\begin{array}{ll}
A_{\mu } = & \mu_{\K}\circ (\mu_{\K}\otimes Id)\circ \lambda _3(f_1\otimes f_2\otimes f_3)
\circ (\Delta \otimes Id)\circ \Delta \\
\smallskip \\
&- \mu_{\K}\circ (Id\otimes \mu_{\K})\circ \lambda _3(f_1\otimes f_2\otimes f_3)\circ 
(Id \otimes \Delta)\circ \Delta. 
\end{array}
$$
and using associativity and commutativity of the multiplication in $\K$, we obtain
$$A_{\mu } =  \mu_{\K}\circ (\mu_{\K}\otimes Id)\circ \lambda _3(f_1\otimes f_2\otimes f_3)
\circ ((\Delta \otimes Id)\circ \Delta-(Id\otimes \Delta )\circ \Delta ).$$
For every $\sigma \in \Sigma _3$ and any $\K$-vector space $E$, we denote by $\Phi _{\sigma }^E$
the action 
$$\Phi _{\sigma }^E :E^{\otimes ^3}\longrightarrow E^{\otimes ^3}$$
defined in the first section :
$$\Phi _{\sigma }^E (e_1\otimes e_2\otimes e_3)=e_{\sigma ^{-1}(1)}\otimes e_{\sigma ^{-1}(2)}
\otimes e_{\sigma ^{-1}(3)}.$$ 
Thus 
$$
\begin{array}{lllllllllllllllllllllll}
\sum_{\sigma \in G_i}(-1)^{\epsilon (\sigma )}A_{\mu }\circ \Phi _{\sigma }^{C^*} &&&&&&&&&&&&&&&&&&&&&&
\end{array}
$$
$$
\begin{array}{ll}
& =  
\mu_{\K}\circ (\mu_{\K}\otimes Id)\circ \lambda _3(f_1\otimes f_2\otimes f_3)
\circ (G_i\circ (\Delta \otimes Id)\circ \Delta-G_i\circ (Id\otimes \Delta )\circ \Delta )\\
\smallskip \\
&=0.
\end{array}$$
$\Box$

\medskip

\noindent {\bf D. The dual space of a finite dimensional $G_i$-associative algebra}

\medskip

\begin{proposition}
The dual vector space of a finite dimensional $G_i$-associative algebra has a $G_i$-associative 
cogebra structure.
\end{proposition}
{\it Proof.} Let $ \mathcal{A} $ be a finite dimensional $G_i$-associative algebra and let 
$\{{e_i},{i=1,...,n}\}$ be a basis of $\mathcal{A} $. If 
$\{f_i\}$ is the dual basis then $\{f_i \otimes f_j\}$ is a basis of 
$\mathcal{A}^* \otimes \mathcal{A}^* $. 
The coproduct $\Delta $ on $\mathcal{A}^*$ is defined by 
$$\Delta (f)=\sum _{i,j} f(\mu (e_i \otimes e_j)) f_i \otimes f_j.$$
In particular
$$\Delta (f_k) = \sum _{i,j} C_{ij}^{k}f_i \otimes f_j$$
where $C_{ij}^k$ are the structure constants of $\mu$ related to the basis $\{{e_i}\}$. 
Then $\Delta$ is the comultiplication of a $G_i$-associative cogebra.

\medskip

\noindent {\bf E. Associated Lie cogebras}

\medskip

Since every $G_i$-associative algebra $(\mathcal{A},\mu )$ is  Lie-admissible, the bilinear map
$\mu (X,Y)-\mu (Y,X)$ is a Lie bracket. Similary we prove a similar  property for $G_i$-associative 
cogebras. Let us recall
the notion of Lie cogebras.  Let $C$ un $\K$-vector space and 
$$\Delta  : C\longrightarrow C\otimes C$$
a linear map satisfying 

1. $\tau \circ \Delta =-\Delta $ where $\tau $ is the permutation $\tau (x\otimes y)=y\otimes x$

2. $\Phi _v \circ  (Id\otimes \Delta )\circ \Delta =0$ where $v= \sum _{\sigma \in G_5}\sigma .$

\noindent Then $(C,\Delta )$ is called a Lie cogebra.

\begin{proposition}
Let $(C,\Delta )$ be a $G_i$-associative cogebra (non necessary counitary). 
Let $\Delta _L$ the linear
map defined by $\Delta _L=\Delta -\tau \circ \Delta$. Then $(C,\Delta _L)$ is a Lie cogebra.
\end{proposition}
{\it Proof.} It is clear that $\tau \circ \Delta_L =-\Delta_L $. 
As $W$ is in the orbit of the vector
$u_i=\sum _{\sigma \in G_i} \sigma $, then any $G_i$-associative cogebra is  
$G_6$-associative cogebra, that is a Lie admissible cogebra. But the identity
$$\Phi _W \circ (Id \otimes \Delta )\circ \Delta - \Phi _W \circ (\Delta \otimes Id )\circ
 \Delta =0$$
is equivalent to
$$
\begin{array}{l}
\Phi _v \circ (Id \otimes \Delta )\circ \Delta - 
\tau \circ \Phi _v \circ (Id \otimes \Delta  )\circ \Delta \\
-\Phi _v \circ  (Id\otimes \tau \circ \Delta )\circ \Delta 
-\tau \circ \Phi _v \circ  (Id \otimes \tau \circ \Delta  )\circ \Delta =0
\end{array}
$$
with $v=Id+c+c^2$, we deduce the result.
 $\Box$

\bigskip

\begin{center}
{\bf IV. THE CONVOLUTION PRODUCT}
\end{center}

\medskip

Let $( \mathcal{A} ,\mu,\eta,G_1)$ be an associative algebra and
$(C, \Delta, \epsilon, G_1)$ an associative cogebra. The convolution product
on the vector space $Hom(C,\mathcal{A})$ is given by
$$
f * g = \mu \circ \lambda_2 \circ  (f\otimes g)\circ \Delta 
$$
for every $f$ and $g$ $\in Hom(C,\mathcal{A})$ 
and $(Hom(C,\mathcal{A}), \star, G_1)$ is an
associative algebra. We want to generalize this property to the other groups $G_i$, $i\geq 2$.
Let us begin by introduce new
classes of associative algebras which appear naturally by studying the duality 
in the corresponding operads (see section V).

\medskip

\noindent {\bf A. The $G_i^!$-algebras and cogebras}

\medskip

\begin{definition}
For $i \geq 2$, a $G_i^!$-algebra is an associative algebra  satisfying :

- for $i=2$ : $x_1.x_2.x_3=x_2.x_1.x_3$

- for $i=3$ : $x_1.x_2.x_3=x_1.x_3.x_2$

- for $i= 4$ : $x_1.x_2.x_3=x_3.x_2.x_1$

- for $i=5$ : $x_1.x_2.x_3=x_2.x_3.x_1=x_3.x_1.x_2$

- for $i=6$ : $x_1.x_2.x_3=x_{\sigma(1)}.x_{\sigma(2)}.x_{\sigma(3)} $ for all $x_1,x_2,x_3$ and $\sigma \in  \Sigma _3$.
\end{definition}
\begin{definition}
For $i \geq 2$, a $G_i^!$-cogebra is a coassociative cogebra, that is a $G_1$-associative 
cogebra,
satisfying :
$$\Phi_{u_i} \circ (Id\otimes \Delta )\circ \Delta =(Id\otimes \Delta )\circ \Delta $$
where $u_i=\sum _{\sigma \in G_i}\sigma^{-1}.$
\end{definition}
\medskip

\noindent {\bf B. Convolution product on $Hom(C,\mathcal{A})$}

\medskip

We will provide $Hom(C,\mathcal{A})$ with a structure of $G_i$-associative algebra.

\begin{proposition}
Let $(\mathcal{A},\mu,\eta ,G_i)$ be a $G_i$-associative algebra and $(C,\Delta,\varepsilon ,G_i^! )$
 a $G_i^!$-cogebra. Then  $(Hom(C,\mathcal{A}), \star,\eta \circ \epsilon, G_i)$ is a
$G_i$-associative algebra where $\star $ is the convolution product :
$$f\star g=\mu \circ \lambda _2(f\otimes g)\circ \Delta .$$
\end{proposition}

\noindent {\it Proof.}  Let us compute the associator $A_*$ of the convolution product.
$$
\begin{array}{ll}
(f_1 \star f_2) \star f_3 & = \mu \ \circ \lambda _2 ((f_1 \star f_2)\otimes f_3 )\,  \circ \,
\Delta
\\ 
& \vspace{0.1cm}\\
& =\mu  \, \circ \, \lambda _2 ((\mu\circ \lambda _2 (f_1 \otimes f_2)\circ  \Delta) \otimes f_3 ) \,  \circ \,
\Delta  \\
& \vspace{0.1cm}\\
& =\mu  \, \circ \,  (\mu \otimes Id)\circ \lambda _3(f_1 \otimes f_2 \otimes
f_3)\circ (\Delta \otimes Id)\, \circ \, \Delta.
\end{array}
$$
Thus 
$$
\begin{array}{ll}
A_{\star }(f_1\otimes f_2\otimes f_3)  &=\mu  \, \circ \,  (\mu \otimes Id)\circ \lambda _3(f_1
 \otimes f_2 \otimes f_3)\circ (\Delta \otimes Id)\, \circ \, \Delta\\
& \vspace{0.1cm}\\
& - \mu  \, \circ \,  (Id \otimes \mu )\circ \lambda _3(f_1 \otimes f_2 \otimes
f_3)\circ (Id\otimes \Delta )\, \circ \, \Delta.
\end{array}
$$
Therefore
$$
\begin{array}{llllllllllllllllllll}
\sum_{\sigma \in G_i}(-1)^{\epsilon (\sigma )}A_{\star } \circ \Phi_{\sigma }^{Hom(C,A)}(f_1\otimes f_2\otimes f_3) 
&&&&&&&&&&&&&&&&&&&
\end{array}
$$
$$\begin{array}{lll}
& =&\mu  \, \circ \,  (\mu \otimes Id)\circ 
(\sum_{\sigma \in G_i}\lambda _3( \Phi_{\sigma }^{Hom(C,A)}(f_1
 \otimes f_2 \otimes f_3)))\circ (\Delta \otimes Id)\, \circ \, \Delta\\
& &\vspace{0.05cm}\\
&& - \mu  \, \circ \,  (Id \otimes \mu )\circ (\sum_{\sigma \in G_i} \lambda _3(
\Phi_{\sigma }^{Hom(C,A)}(f_1 \otimes f_2 \otimes
f_3)))\circ (Id\otimes \Delta )\, \circ \, \Delta.
\end{array}
$$
But
$$\lambda _3(
\Phi_{\sigma }^{Hom(C,A)}(f_1 \otimes f_2 \otimes
f_3))=\Phi_{\sigma }^{A}\circ \lambda _3(f_1 \otimes f_2 \otimes
f_3)\circ \Phi_{\sigma^{-1} }^{C}.$$
This gives 
$$
\begin{array}{llllllllllllllllllll}
\sum_{\sigma \in G_i}(-1)^{\epsilon (\sigma )}A_{\star } \circ \Phi_{\sigma }^{Hom(C,A)}(f_1\otimes f_2\otimes f_3) 
&&&&&&&&&&&&&&&&&&&
\end{array}
$$
$$\begin{array}{lll}
& =&\mu  \, \circ \,  (\mu \otimes Id)\circ 
(\sum_{\sigma \in G_i}
\Phi_{\sigma }^{A}\circ \lambda _3(f_1 \otimes f_2 \otimes
f_3))\circ \Phi_{\sigma^{-1} }^{C}\circ 
(\Delta \otimes Id)\, \circ \, \Delta\\
&&\vspace{0.05cm}\\
&& - \mu  \, \circ \,  (Id \otimes \mu )\circ (\sum_{\sigma \in G_i} 
\Phi_{\sigma }^{A}\circ \lambda _3(f_1 \otimes f_2 \otimes
f_3))\circ \Phi_{\sigma^{-1} }^{C}\circ (Id\otimes \Delta )\, \circ \, \Delta.
\end{array}
$$
\smallskip

\noindent As $\Delta $ is coassociative 
$$(\Delta \otimes Id)\, \circ \, \Delta=(Id\otimes \Delta )\, \circ \, \Delta$$
and the $G_i^!$-cogebra structure implies
$$\Phi_{u_i }^{C}\circ (Id\otimes \Delta )\, \circ \, \Delta=
\Phi_{u_i }^{C}\circ 
(\Delta \otimes Id)\, \circ \, \Delta=
(\Delta \otimes Id)\, \circ \, \Delta.$$
Then
$$
\begin{array}{llllllllllllllllllll}
\sum_{\sigma \in G_i}(-1)^{\epsilon (\sigma )}A_{\star } \circ \Phi_{\sigma }^{Hom(C,A)}(f_1\otimes f_2\otimes f_3) 
&&&&&&&&&&&&&&&&&&&
\end{array}
$$
$$
\begin{array}{lll}
& =&\mu  \, \circ \,  (\mu \otimes Id)\circ 
(\sum_{\sigma \in G_i}
\Phi_{\sigma }^{A}\circ \lambda _3(f_1 \otimes f_2 \otimes
f_3))\circ (\Delta \otimes Id)\, \circ \, \Delta\\
&&\vspace{0.05cm}\\
&& - \mu  \, \circ \,  (Id \otimes \mu )\circ (\sum_{\sigma \in G_i} 
\Phi_{\sigma }^{A}\circ \lambda _3(f_1 \otimes f_2 \otimes
f_3))\circ (\Delta \otimes Id)\, \circ \, \Delta\\
& &\smallskip\\
&=&\sum_{\sigma \in G_i}A_{\mu }\circ \Phi_{\sigma }^{A}\circ \lambda _3(f_1 \otimes f_2 \otimes
f_3)\circ (\Delta \otimes Id)\, \circ \, \Delta \\
&&\smallskip \\
&=&0.
\end{array}
$$
\smallskip

\noindent This proves the proposition. $\Box$

\bigskip

\begin{center}
{\bf V. TENSOR PRODUCT of  $\Sigma _3$-ASSOCIATIVE ALGEBRAS }
\end{center}

\medskip

\noindent {\bf A. Non monoidal categories}
We know that the tensor product of associative algebras can be provided with an associative algebra structure. In other words, the category
of associative algebras is monoidal and closed for the tensor product. This is 
generally not true for other categories of $\Sigma _3$-associative
algebras. 

\begin{proposition}
Let $(\mathcal{A},\mu _A)$ and $(\mathcal{B},\mu _B)$ be two $\Sigma_3$-associative algebras 
respectively defined by the relations $A_{\mu_A} \circ \Phi_v =0$ and 
$A_{\mu_B} \circ \Phi_w =0$. Then $(\mathcal{A}\otimes_{\K} \mathcal{B},\mu_A \otimes \mu_B)$
is a $\Sigma_3$-associative algebra if and only if  $\mathcal{A}$ and $\mathcal{B}$ are
associative algebras (i.e $G_1$-associative algebras).
\end{proposition}

\noindent {\it Proof.} Let $A_\mu$ be the associator of the law $\mu$ and
$A_{\mu}^L (x_1\otimes x_2 \otimes x_3)=\mu(\mu(x_1, x_2),x_3)$ and 
$A_{\mu}^R =A_{\mu}^L -A_{\mu}$. By hypothesis, $\mathcal{A}$ (resp. $\mathcal{B}$)
is defined by $A_{\mu_A} \circ \Phi_v=0$ 
(resp. $A_{\mu_B} \circ \Phi_w=0)$. Let us suppose that $\mathcal{A} \otimes \mathcal{B}$ is a $\Sigma_3$-associative algebra. 
There exists
$v' \in \mathbb{K}$ such that $A_{\mu_{A \otimes B}} \circ \Phi_{v'}=0.$ By taking $v'=\sum_{i=1}^6 \gamma_i \sigma_i$
the last condition becomes
\begin{eqnarray}
\label{3} \sum_{i=1}^6 \gamma_i [A_{\mu_{A \otimes B}} \circ \Phi_{\sigma_i}]=0
\end{eqnarray}
and 
\begin{eqnarray}
\label{4} \sum_{i=1}^6 \gamma_i [A_{\mu_{A}}^L \circ \Phi_{\sigma_i} \otimes A_{\mu_{B}}^L\circ \Phi_{\sigma_i}
-A_{\mu_{A}}^R\circ \Phi_{\sigma_i} \otimes A_{\mu_{B}}^R\circ \Phi_{\sigma_i}]=0.
\end{eqnarray}
Let us denote by
$e_i=A_{\mu_{A}}^L\circ \Phi_{\sigma_i},\widetilde{e_i}=A_{\mu_{A}}^R\circ \Phi_{\sigma_i},f_i=A_{\mu_{B}}^L\circ \Phi_{\sigma_i}$ 
and $\widetilde{f_i}=A_{\mu_{B}}^R\circ \Phi_{\sigma_i}.$ The vectors $e_i$ and 
$\widetilde{e_i}$ belong to $Hom(\mathcal{A} ^{\otimes ^3}, \mathcal{A})$ and the vectors 
$f_i$ and $\widetilde{f_i}$ to $Hom(\mathcal{B} ^{\otimes ^3}, \mathcal{B})$.  
Equation \ref{4} becomes :
\begin{eqnarray}
\label{5} \quad \sum_{i=1}^6 \gamma_i [e_i \otimes f_i - \widetilde{e_i} \otimes \widetilde{f_i}]=0.
\end{eqnarray}
From the definition of the algebras $\mathcal{A}$ and $\mathcal{B}$ we have:
$$\sum_{i=1}^6 a_i (e_i -\widetilde{e_i})=0 \ {\rm and} \ \sum_{i=1}^6 b_i (f_i -\widetilde{f_i})=0$$
if $v=\sum_{i=1}^6 a_i \sigma_i$ and $w=\sum_{i=1}^6 b_i \sigma_i.$ Suppose that $dim F_w=k$ with $0\leq k \leq 6$. Then the rank of the vectors $\left\{ f_i,\widetilde{f_i} \right\}$ is equal to $(6+k).$ We can suppose that
$\left\{ f_1,...,f_6,\widetilde{f_1},...,\widetilde{f_k} \right\}$ are independent. This implies:
$$
\left\{
\begin{array}{l}
\widetilde{f}_{k+1}=\rho_1^{k+1}f_1+...+\rho_6^{k+1}f_6+\widetilde{\rho}_1^{k+1}\widetilde{f}_1+...+\widetilde{\rho}_k^{k+1}\widetilde{f}_k \\
: \\
: \\
\widetilde{f}_{6}=\rho_1^{6}f_1+...+\rho_6^{6}f_6+\widetilde{\rho}_1^{6}\widetilde{f}_1+...+\widetilde{\rho}_k^{6}
\widetilde{f}_k. \\
\end{array}
\right.
$$
The equation \ref{5} becomes:
$$
\sum_{i=1}^6 e_i' \otimes f_i +\sum_{i=1}^k e_i^{\prime \prime} \otimes \widetilde{f}_k=0
$$
and the independence of the vectors $\left\{ f_1,...,f_6,\widetilde{f}_1,...,\widetilde{f}_k \right\}$ implies
$e'_1=...=e'_6=e^{\prime\prime}_1=...=e^{\prime\prime}_k=0.$ Then we deduce
$$
\left\{
\begin{array}{l}
\gamma_1 e_1- \gamma_{k+1}\rho_1^{k+1}\widetilde{e}_1-...-\gamma_{6}\rho_1^{6}\widetilde{e}_6=0 \\
: \\
: \\
\gamma_6 e_6- \gamma_{k+1}\rho_6^{k+1}\widetilde{e}_1-...-\gamma_{6}\rho_6^{6}\widetilde{e}_6=0 \\
\gamma_1 \widetilde{e_1}+ \gamma_{k+1}\widetilde{\rho}_1^{k+1}\widetilde{e}_{k+1}+...+\gamma_{6}\widetilde{\rho}_1^{6}\widetilde{e}_6=0 \\
: \\
: \\
\gamma_k \widetilde{e}_k+ \gamma_{k+1}\widetilde{\rho}_k^{k+1}\widetilde{e}_{k+1}+...+\gamma_{6}\widetilde{\rho}_k^{6}\widetilde{e}_6=0. \\
\end{array}
\right.
$$
Let us suppose that $k\neq 0.$ Then the conditions on the vector $v$ are of the type
$\sum_{i=1}^6 a_i(e_i -\widetilde{e}_i)=0$ and we have succesively $\gamma_1=...=\gamma_k=0$ and $\gamma_{k+1}\widetilde{\rho}_i^{k+1}=...=\gamma_{6}\widetilde{\rho}_i^{6}=0$ for $i=1,...,6$. If one of $\gamma_j$ for $j=k+1,...,6$ is not 0 then $\widetilde{\rho}_j^{i}=0$ for $i=1,...,6$. Thus
$$\widetilde{f}_j=\rho_1^{j}f_1+...+\rho_6^{j}f_6$$
and, as for every $i$ there exists $\sigma \in \Sigma _3$ such that $\widetilde{f}_i=\widetilde{f}_j \circ \Phi _{\sigma}$, we have also
$$\widetilde{f}_i=\rho_1^{i}f_1+...+\rho_6^{i}f_6$$
for every $i=1,..,6$. This implies $k=0$ which is impossible from the hypothesis.
Thus $k=0$ and the vectors $\left\{ f_i, \widetilde{f}_i \right\}_{i=1,...,6}$ is of rank 6.
The only possible relations are then $f_i=\widetilde{f}_i$ and $\mathcal{B}$ is associative.
We deduce the associativity of $\mathcal{A}$.
 $\Box$

\medskip

\noindent {\bf B. $G_i!$-algebras}

\medskip

In \cite{G.R}, we have studied the operads corresponding to the $G_i$-associative algebras. More precisely,
the  operad  $G_{i}-\mathcal{A}ss$ is the quadratic operad defined by the relations
$$\sum_{\sigma  \in G_i}(-1)^{\epsilon (\sigma )}
(x_{\sigma ^{-1}(1)}.x_{\sigma ^{-1}(2)}).x_{\sigma ^{-1}(3)}-
x_{\sigma ^{-1}(1)}.(x_{\sigma ^{-1}(2)}).x_{\sigma ^{-1}(3)}).$$
Then an algebra on the operad $G_{i}-\mathcal{A}ss$ is a $G_i$-associative algebra.
In \cite{G.R}, we have also computed the dual operads $(G_{i}-\mathcal{A}ss)!$ and proved
that only $G_{1}-\mathcal{A}ss$ (classicaly noted $\mathcal{A}ss$ is autodual that is 
$(G_{i}-\mathcal{A}ss)=(G_{i}-\mathcal{A}ss)!$. We call a $G_i!$-algebra an algebra on
the quadratic operad $(G_{i}-\mathcal{A}ss)!$. In definition 13 we have given their structures.

\medskip

\noindent {\bf C. Tensor products}

\medskip

\begin{theorem}
If $\mathcal{A}$ is a $G_i$-assocative algebra and $\mathcal{B}$ a $G_i^!$-algebra 
(with the same indice) then $\mathcal{A} \otimes \mathcal{B}$ can be provided with a $G_i$-algebra structure
for $i=1,...,6$. 

\end{theorem}
{\it Proof.} Let us consider on  $\mathcal{A} \otimes \mathcal{B}$  the classical tensor
 product
$$\mu _A\otimes \mu _B((a_1 \otimes b_1)\otimes (a_2 \otimes b_2))
=\mu _A(a_1 \otimes a_2) \otimes \mu _B(b_1 \otimes b_2).$$
To simplify, we denote by $\mu $ the product $\mu _A\otimes \mu _B$.
As $\mathcal{B}$  is an associative algebra, the associator $A_{\mu }$ satisfies:
$$A_{\mu }((a_1 \otimes b_1)\otimes (a_2 \otimes b_2)\otimes (a_3 \otimes b_3))=
A_{\mu _A}(a_1\otimes a_2\otimes a_3) \otimes 
\mu _B\circ (\mu _B \otimes Id)(b_1 \otimes b_2 \otimes b_3).$$
Therefore
$$
\begin{array}{lllllllllll}
\sum _{\sigma \in G_i}(-1)^{\epsilon (\sigma )} A_{\mu }\circ \Phi_{\sigma }
((a_1 \otimes b_1)\otimes (a_2 \otimes b_2)\otimes (a_3 \otimes b_3)) & &&&&&&&&&
\end{array}
$$
$$
\begin{array}{ll}
&=\sum _{\sigma \in G_i}(-1)^{\epsilon (\sigma )}A_{\mu_A }\circ \Phi_{\sigma }
(a_1\otimes a_2\otimes a_3) \otimes 
\mu _B\circ (\mu _B \otimes Id) \circ \Phi_{\sigma} 
(b_1 \otimes b_2 \otimes b_3).
\end{array}
$$
But $\mathcal{B}$ a $G_i^!$-algebra. Then  
$$\mu _B\circ (\mu _B \otimes Id) \circ \Phi_{\sigma} 
(b_1 \otimes b_2 \otimes b_3)=\mu _B\circ (\mu _B \otimes Id) 
(b_1 \otimes b_2 \otimes b_3).$$
We obtain
$$
\begin{array}{lllllllllll}
\sum _{\sigma \in G_i}(-1)^{\epsilon (\sigma )} A_{\mu }\circ \Phi_{\sigma }
((a_1 \otimes b_1)\otimes (a_2 \otimes b_2)\otimes (a_3 \otimes b_3)) & &&&&&&&&&
\end{array}
$$
$$
\begin{array}{ll}
&=(\sum _{\sigma \in G_i}(-1)^{\epsilon (\sigma )}A_{\mu_A }\circ \Phi_{\sigma }
(a_1\otimes a_2\otimes a_3)) \otimes 
\mu _B\circ (\mu _B \otimes Id)  
(b_1 \otimes b_2 \otimes b_3) \\
\medskip
& = 0.
\end{array}
$$
$\Box$

\medskip

\noindent{\bf D. The categories $G_i-ASS$ }

\medskip

Let $G_i-ASS$ and  $(G_i -ASS)^{!}$ the categories whose objects are respectively
 $G_i$-associative algebras and $G_i^!$-algebras and the morphismes are
the homomorphismes of algebras. The previous theorem can be interpreted as follows: 
Let  $\mathcal{A}$ be a $G_i$-associative algebra. Then  $\mathcal{A} \otimes -$ is a covariant functor
$$ \mathcal{A} \otimes - : (G_i -ASS)^{!} \longrightarrow G_i -ASS. $$

\medskip

\noindent{\bf Thanks}: I thank the referee to point up a mistake on the section IV in the first version
and to propose a correction.


\begin{thebibliography}{99}
\bibitem{Al}  Albert A.A.,  Power-associative rings. {\it Trans. Amer. Math. Soc.} {\bf 64} (1948), 552-593.

\bibitem{B.H} Bremmer M., Hentzel I., Identities for the associator in alternative algebras. Preprint.


\bibitem{Dz} Dzhumadil'daev A., Novikov-Jordan algebras. {\it Comm. in Algebra} {\bf 11}, (2002), 5207-5240.

\bibitem{G} Gerstenhaber M., The cohomology structure of an associative ring,
{\it Ann of math.} {\bf 78}, 2,  (1963), 267-288.

\bibitem{G.K} Ginzburg V., Kapranov M., Koszul duality for operads. {\it
Duke Math Journal}. {\bf 76} 1, (1994), 203-272.

\bibitem{E} Elduque A., Myung H.C.,  {\it Mutations of Alternative Algebras.} Kluwer Academic Press. Math and its applications 278. 1994.

\bibitem{G.R} Goze M., Remm E., Lie-admissible algebras and operads. {\it Journal of Algebra}
{\bf 273/1} (2004), 129-152.

\bibitem{G.R 2} Goze M., Remm E., A class of nonassociative algebras. Preprint Mulhouse 2003 and xxx RA/0309015


\bibitem{N} Nijenhuis A., Sur une classe de propri\'et\'es communes
\`{a} quelques types diff\'erents d'alg\`ebres. {\it Enseignement Math.} {\bf
14} (2), (1970), 225--277.

\bibitem{R} Remm E., Op\'erades Lie-admissibles.
{\it C. R. Math. Acad. Sci. Paris}  {\bf 334}  no. 12, (2002), 1047--1050.

\bibitem{S} Schafer R.D. {\it Introduction to Nonassociative Algebras}. Academic Press. 1966
\end{thebibliography}
\end{document}